\documentclass{amsart}
\usepackage{latexsym,amsmath,amssymb}

\theoremstyle{definition}
\theoremstyle{remark}

\numberwithin{equation}{section}

\begin{document}

\title[Weighted conditional type operators ]
{Centered weighted conditional type operators }

\author{\sc\bf Y. Estaremi  }
\address{\sc Y. Estaremi  }
\email{yestaremi@pnu.ac.ir}

\address{Department of Mathematics, Payame Noor University, p. o. box: 19395-3697, Tehran, Iran.}

\thanks{}

\thanks{}



\date{}

\dedicatory{}

\commby{}

\begin{abstract}
In this paper, we give some necessary and sufficient conditions
for weighted conditional expectation type operators on
$L^{2}(\Sigma)$ to be centered. Also, we investigate the relation
between normal and centered weighted conditional type operators.
Finally we give some applications.

 \noindent {}
\end{abstract}

\maketitle

\section{ \sc\bf Introduction and Preliminaries}

Operators in function spaces defined by conditional expectations
were first studied, among others, by S - T.C. Moy \cite{mo}, Z.
Sidak \cite{s} and H.D. Brunk \cite{b} in the setting of $L^p$
spaces. Conditional expectation operators on various function
spaces exhibit a number of remarkable properties related to the
underlying structure of the given function space or to the metric
structure when the function space is equipped with a norm. P.G.
Dodds, C.B. Huijsmans and B. de Pagter \cite{dhd} linked these
operators to averaging operators defined on abstract spaces
earlier by J.L. Kelley \cite{k}, while A. Lambert \cite{l} studied
their link to classes of multiplication operators which form
Hilbert $C^*$-modules. J.J. Grobler and B. de Pagter \cite{gd}
showed that the classes of partial integral operators, studied by
A.S. Kalitvin and others \cite{afkz, akn, akz, dou, kz}, were a
special case of conditional expectation operators. Recently, J.
Herron studied operators $EM_u$ on $L^p$
spaces in \cite{ her}. \\

Let $f\in L^0(\Sigma)$, then $f$ is said to be conditionable with
respect to $E$ if $f\in\mathcal{D}(E):=\{g\in L^0(\Sigma):
E(|g|)\in L^0(\mathcal{A})\}$. Throughout this paper we take $u$
and $w$ in $\mathcal{D}(E)$. Let $(X,\Sigma,\mu)$ be a complete
$\sigma$-finite measure space. For any sub-$\sigma$-finite algebra
$\mathcal{A}\subseteq
 \Sigma$, the $L^2$-space
$L^2(X,\mathcal{A},\mu_{\mid_{\mathcal{A}}})$ is abbreviated  by
$L^2(\mathcal{A})$, and its norm is denoted by $\|.\|_2$. All
comparisons between two functions or two sets are to be
interpreted as holding up to a $\mu$-null set. The support of a
measurable function $f$ is defined as $S(f)=\{x\in X; f(x)\neq
0\}$. We denote the vector space of all (equivalence classes of)
almost everywhere finite valued measurable functions on $X$ by
$L^0(\Sigma)$.

\vspace*{0.3cm} For a sub-$\sigma$-finite algebra
$\mathcal{A}\subseteq\Sigma$, the conditional expectation operator
associated with $\mathcal{A}$ is the mapping $f\rightarrow
E^{\mathcal{A}}f$, defined for all non-negative, measurable
function $f$ as well as for all $f\in L^2(\Sigma)$, where
$E^{\mathcal{A}}f$, by the Radon-Nikodym theorem, is the unique
$\mathcal{A}$-measurable function satisfying
$$\int_{A}fd\mu=\int_{A}E^{\mathcal{A}}fd\mu, \ \ \ (A\in \mathcal{A}).$$
As an operator on $L^{2}({\Sigma})$, $E^{\mathcal{A}}$ is
idempotent and $E^{\mathcal{A}}(L^2(\Sigma))=L^2(\mathcal{A})$. If
there is no possibility of confusion, we write $E(f)$ in place of
$E^{\mathcal{A}}(f)$.

\vspace*{0.2cm}\noindent A detailed discussion of the properties
of these operators may be found in \cite{rao}. We recall that an
$\mathcal{A}$-atom of the measure $\mu$ is an element
$A\in\mathcal{A}$ with $\mu(A)>0$ such that for each
$F\in\mathcal{A}$, if $F\subseteq A$, then either $\mu(F)=0$ or
$\mu(F)=\mu(A)$. A measure space $(X,\Sigma,\mu)$ with no atoms is
called a non-atomic measure space. It is well-known fact that
every $\sigma$-finite measure space $(X,
\mathcal{A},\mu_{\mid_{\mathcal{A}}})$ can be partitioned uniquely
as $X=\left (\bigcup_{n\in\mathbb{N}}A_n\right )\cup B$, where
$\{A_n\}_{n\in\mathbb{N}}$ is a countable collection of pairwise
disjoint $\mathcal{A}$-atoms and $B$, being disjoint from each
$A_n$, is non-atomic (see \cite{z}).\\

An operator $A$ on a Hilbert space is {\it centered} if the family
of operators $\{A^{\ast^n}A^n, A^kA^{\ast^k}: n,k\geq0\}$ is
commutative \cite{mm}.\\

In \cite{e, ej} we investigated some classic properties of
multiplication conditional expectation operators $M_wEM_u$ on
$L^p$-spaces. In this paper we will be concerned with
characterizing weighted conditional expectation type operators and
their  Aluthage transformations on $L^2(\Sigma)$ in terms of
membership in the class of centered operators and the relation
between normal and centered weighted conditional type operators.
Finally, we get some necessary and sufficient conditions for integral operators to be centered and some applications in their spectra.\\

\section{ \sc\bf Centered and normal weighted conditional type operators}

In the first we reminisce some theorems that we have proved in
\cite{ej}.

\vspace*{0.3cm} {\bf Theorem 2.1.} The operator $T=M_wEM_u$ is
bounded on $L^{2}(\Sigma)$ if and only if
$(E|w|^{2})^{\frac{1}{2}}(E|u|^{2})^{\frac{1}{2}} \in
L^{\infty}(\mathcal{A})$, and in this case its norm is given by
$\|T\|=\|(E(|w|^{2}))^{\frac{1}{2}}(E(|u|^{2}))^{\frac{1}{2}}\|_{\infty}$.\\

\vspace*{0.3cm} {\bf Lemma 2.2.} Let $T=M_wEM_u$ be a bounded
operator on $L^{2}(\Sigma)$ and let $p\in (0,\infty)$. Then
$$(T^{\ast}T)^{p}=M_{\bar{u}(E(|u|^{2}))^{p-1}\chi_{S}(E(|w|^{2}))^{p}}EM_{u}$$
and
$$(TT^{\ast})^{p}=M_{w(E(|w|^{2}))^{p-1}\chi_{G}(E(|u|^{2}))^{p}}EM_{\bar{w}},$$
where $S=S(E(|u|^2))$ and $G=S(E(|w|^2))$.

\vspace*{0.3cm} {\bf Theorem 2.3.} The unique polar decomposition
of bounded operator $T=M_wEM_u$ is $U|T|$, where

$$|T|(f)=\left(\frac{E(|w|^{2})}{E(|u|^{2})}\right)^{\frac{1}{2}}\chi_{S}\bar{u}E(uf)$$
and
 $$U(f)=\left(\frac{\chi_{S\cap
 G}}{E(|w|^{2})E(|u|^{2})}\right)^{\frac{1}{2}}wE(uf),$$
for all $f\in L^{2}(\Sigma)$.

\vspace*{0.3cm} {\bf Theorem 2.4.} The Aluthge transformation of
$T=M_wEM_u$ is
$$\widehat{T}(f)=\frac{\chi_{S}E(uw)}{E(|u|^{2})}\bar{u}E(uf), \ \ \ \ \ \  \ \ \ \  \ \  \ f\in L^{2}(\Sigma).$$
\\
 It is shown in \cite{mm} that if $A$ is an
operator such that, for each positive $n$, $A^n$ has polar
decomposition $U_nP_n$, then $A$ is centered if and only if for
each positive $n$, $U_n=U^n_1$. In the sequel we will characterize
centered weighted conditional type operators.

\vspace*{0.3cm} {\bf Theorem 2.5.} Consider the weighted
conditional type operator $M_wEM_u:L^2(\Sigma)\rightarrow
L^2(\Sigma)$. Then\\

(a) If $M_wEM_u$ is centered, then $|E(uw)|^2=E(|u|^2)E(|w|^2)$ on $S(E(uw)E(w)E(u))$.\\

(b) If $|E(uw)|^2=E(|u|^2)E(|w|^2)$, then $M_wEM_u$ is centered.\\

\vspace*{0.3cm}{\bf Proof.} (a) By induction we have
$$T^nf=(E(uw))^{n-1}wE(uf),  \ \ \ \ f\in L^2(\Sigma), \ \  n\in
\mathbb{N}.$$  Now, by Theorem 2.3 we obtain

$$U_n(f)=\frac{\chi_{H}
E(uw)^{n-1}}{(E(|u|^2))^{\frac{1}{2}(E(|w|^2))^{\frac{1}{2}}}
|E(uw)|^{n-1}}wE(uf), \ \ \ U_1(f)=\frac{\chi_{S\cap
 G}}{(E(|u|^{2}))^{\frac{1}{2}}(E(|w|^{2}))^{\frac{1}{2}}}wE(uf)$$

 and
 $$U^n_1(f)=\frac{\chi_{S\cap G}E(uw)^{n-1}}{(E(|u|^{2}))^{\frac{n}{2}}(E(|w|^{2}))^{\frac{n}{2}}}wE(uf)$$
 for all $f\in L^2(\Sigma)$ and $H=S(E(uw))$. This implies that if $U_n=U^n_1$,
 then
$$\frac{\chi_{H}
E(uw)^{n-1}}{(E(|u|^2))^{\frac{1}{2}(E(|w|^2))^{\frac{1}{2}}}
|E(uw)|^{n-1}}wE(uf)=\frac{\chi_{S\cap
G}E(uw)^{n-1}}{(E(|u|^{2}))^{\frac{n}{2}}(E(|w|^{2}))^{\frac{n}{2}}}wE(uf)$$
for all $f\in L^2(\Sigma)$. So, for positive element $a\in
L^2(\mathcal{A})$,

$$\frac{\chi_{H}
E(uw)^{n-1}}{(E(|u|^2))^{\frac{1}{2}(E(|w|^2))^{\frac{1}{2}}}
|E(uw)|^{n-1}}wE(ua)=\frac{\chi_{S\cap
G}E(uw)^{n-1}}{(E(|u|^{2}))^{\frac{n}{2}}(E(|w|^{2}))^{\frac{n}{2}}}wE(ua)$$
and so

$$(\frac{\chi_{H}
E(uw)^{n-1}}{(E(|u|^2))^{\frac{1}{2}(E(|w|^2))^{\frac{1}{2}}}
|E(uw)|^{n-1}}-\frac{\chi_{S\cap
G}E(uw)^{n-1}}{(E(|u|^{2}))^{\frac{n}{2}}(E(|w|^{2}))^{\frac{n}{2}}})wE(u)a=0.$$

Then

$$
|E(uw)|^{n-1}(E(|u|^2))^{\frac{1}{2}}(E(|w|^2))^{\frac{1}{2}}=
(E(|u|^{2}))^{\frac{n}{2}}(E(|w|^{2}))^{\frac{n}{2}}$$

on $S(E(w)E(u))\cap H$. Thus

  $|E(uw)|^2=E(|u|^2)E(|w|^2)$ on $S(E(w)E(u))\cap H$.\\

(b) Suppose that $|E(uw)|^2=E(|u|^2)E(|w|^2)$. By part (a) and
direct computation shows that $U_n=U^n_1$. Thus $M_wEM_u$ is
centered.\\

\vspace*{0.3cm} {\bf Corollary 2.6.} If $S(E(u)E(w))=S\cap G=H$,
then the operator $M_wEM_u$ on $L^2(\Sigma)$ is centered if and
only if
$|E(uw)|^2=E(|u|^2)E(|w|^2)$.\\

\vspace*{0.3cm} {\bf Corollary 2.7.} Consider the weighted
conditional type operator $EM_u:L^2(\Sigma)\rightarrow
L^2(\Sigma)$. Then\\

(1) If $EM_u$ is centered, then $|E(u)|^2=E(|u|^2)$ on $S(E(u))$.\\

(2) If $|E(u)|^2=E(|u|^2)$, then $EM_u$ is centered.

\vspace*{0.3cm} {\bf Corollary 2.8.} If $S(E(u))=S(E(|u|^2))$,
then the operator $EM_u$ on $L^2(\Sigma)$ is centered if and only
if
$u\in L^0(\mathcal{A})$.\\

Recall that each operator $A$ on a Hilbert space $\mathcal{H}$ is
called normal if $A^{\ast}A=AA^{\ast}$. In the sequel some
necessary and sufficient conditions for normality will be presented.\\

\vspace*{0.3cm} {\bf Theorem 2.9.} Let $T=M_wEM_u$ be a bounded
operator on $L^{2}(\Sigma)$, then

(a) If
$(E(|u|^2))^{\frac{1}{2}}\bar{w}=u(E(|w|^2))^{\frac{1}{2}},$ then
$T$ is normal.

(b) If $T$ is normal, then $|E(u)|^2E(|w|^2)=|E(w)|^2E(|u|^2)$.\\

{\bf Proof.} (a) Applying lemma 2.2 we have
$$T^{\ast}T-TT^{\ast}=M_{\bar{u}E(|w|^{2})}EM_{u}-M_{wE(|u|^{2})}EM_{\bar{w}}.$$
So for every $f\in L^{2}(\Sigma)$,
$$\langle T^{\ast}T-TT^{\ast}(f),f\rangle=$$
$$\int_{X}E(|w|^2)E(uf)\bar{uf}-E(|u|^2)E(\bar{w}f)w\bar{f}d\mu$$
$$=\int_{X}|E(u(E(|w|^2))^{\frac{1}{2}}f)|^2-|E((E(|u|^2))^{\frac{1}{2}}\bar{w}f)|^2d\mu.$$
This implies that if
$$(E(|u|^2))^{\frac{1}{2}}\bar{w}=u(E(|w|^2))^{\frac{1}{2}},$$
then for all $f\in L^2(\Sigma)$, $\langle
T^{\ast}T-TT^{\ast}(f),f\rangle=0$, thus $T^{\ast}T=TT^{\ast}$.\\

(b) Suppose that $T$ is normal. By (a), for all $f\in L^2(\Sigma)$
we have
$$\int_{X}|E(u(E(|w|^2))^{\frac{1}{2}}f)|^2-|E((E(|u|^2))^{\frac{1}{2}}\bar{w}f)|^2d\mu=0.$$
Let $A\in \mathcal{A}$, with $0<\mu(A)<\infty$. By replacing $f$
to $\chi_{A}$, we have

$$\int_{A}|E(u(E(|w|^2))^{\frac{1}{2}})|^2-|E((E(|u|^2))^{\frac{1}{2}}\bar{w})|^2d\mu=0$$
and so
$$\int_{A}|E(u)|^2E(|w|^2)-|E(w)|^2E(|u|^2)d\mu=0.$$
Since $A\in \mathcal{A}$ is arbitrary, then
$|E(u)|^2E(|w|^2)=|E(w)|^2E(|u|^2)$.$\Box$

\vspace*{0.3cm} {\bf Corollary 2.10.} The operator $EM_u$ on
$L^2(\Sigma)$ is normal if and only if $u\in
L^{\infty}(\mathcal{A})$.\\

\vspace*{0.3cm} {\bf Corollary 2.11.} Consider the weighted
conditional type operator $T=EM_u:L^2(\Sigma)\rightarrow
L^2(\Sigma)$. If $S(E(u))=S(E(|u|^2))$, then the following
conditions are equivalent:\\

(1) $T$ is centered.\\

(2) $T$ is normal.\\

(3) $u \in L^{\infty}(\mathcal{A})$.\\

\vspace*{0.3cm} {\bf Example 2.12.} Let $X=[0,1]\times [0,1]$,
$d\mu=dxdy$, $\Sigma$  the  Lebesgue subsets of $X$ and let
$\mathcal{A}=\{A\times [0,1]: A \ \mbox{is a Lebesgue set in} \
[0,1]\}$. Then, for each $f$ in $L^2(\Sigma)$, $(Ef)(x,
y)=\int_0^1f(x,t)dt$, which is independent of the second
coordinate. Now, if we take $u(x,y)=y^{\frac{x}{8}}$ and $w(x,
y)=\sqrt{(4+x)y}$, then $E(|u|^2)(x,y)=\frac{4}{4+x}$ and
$E(|w|^2)(x,y)=\frac{4+x}{2}$. So, $E(|u|^2)(x,y)E(|w|^2)(x,y)=2$
and $|E(uw)|^2(x,y)=64\frac{4+x}{(x+12)^2}$. Direct computations
shows that
$$E(|u|^2)(x,y)E(|w|^2)(x,y)\leq|E(uw)|^2(x,y).$$

 Since $T=M_wEM_u$ is bounded then we have
$$E(|u|^2)(x,y)E(|w|^2)(x,y)\geq |E(uw)|^2(x,y),$$ this implies that $E(|u|^2)(x,y)E(|w|^2)(x,y)=
|E(uw)|^2(x,y)$. So by Theorem 2.5 the operator $T=M_wEM_u$ is
centered.\\

\vspace*{0.3cm} {\bf Example 2.13.} Let $X=[-1,1]$,
$d\mu=\frac{1}{2}dx$ and $\mathcal{A}=\langle\{(-a,a):0\leq
a\leq1\}\rangle$ ($\sigma$-algebra generated by symmetric
intervals). Then
 $$E^{\mathcal{A}}(f)(x)=\frac{f(x)+f(-x)}{2}, \ \ x\in X,$$
 whenever $E^{\mathcal{A}}(f)$ is defined. Let $u(x)=e^x$, then
 $E(u)(x)=\cosh(x)$ and $E(|u|^2)(x)=\cosh(2x)$. It is clear that
 $S(E(u))=X$. So, $u$ is not $\mathcal{A}$-measurable. Thus by
 Corollary 2.7 we conclude that the operator $T=EM_u$ is not
 centered. If $T=EM_u$ is centered, then $u$ should be
 $\mathcal{A}$-measurable, but $u$ is not $\mathcal{A}$-measurable. \\

\section{ \sc\bf Applications}

In this section, we shall denote by $\sigma_{p}(T)$,
$\sigma_{jp}(T)$, the point spectrum
of $T$, the joint point spectrum of $T$, respectively. A complex number $\lambda\in \mathbb{C}$ is said
to be in the point spectrum $\sigma_{p}(T)$ of the operator $T$,
if there is a unit vector $x$ satisfying $(T-\lambda)x=0$. If in
addition, $(T^{\ast}-\bar{\lambda})x=0$, then $\lambda$ is said to
be in the joint spectrum $\sigma_{jp}(T)$ of $T$.\\

If $A, B\in \mathcal{B}(\mathcal{H})$, then it is well known that
$\sigma_{p}(AB)\setminus\{0\}=\sigma_{p}(BA)\setminus\{0\}$
and
$\sigma_{jp}(AB)\setminus\{0\}=\sigma(BA)_{jp}\setminus\{0\}$.\\

Let
$A_{\lambda}=\{x\in X:E(u)(x)=\lambda\}$, for $\lambda\in
\mathbb{C}$. Suppose that $\mu(A_{\lambda})>0$. Since
$\mathcal{A}$ is $\sigma$-finite,  there exists an
$\mathcal{A}$-measurable subset $B$ of $A_{\lambda}$ such that
$0<\mu(B)<\infty$, and $f=\chi_{B}\in L^p(\mathcal{A})\subseteq
L^p(\Sigma)$. Now
$$EM_u(f)-\lambda f=E(u)\chi_{B}-\lambda \chi_{B}=0.$$ This
implies that $\lambda\in \sigma_{p}(EM_u)$.\\

If there exists $f\in L^p(\Sigma)$ such that $f\chi_{C}\neq 0$
$\mu$-a.e, for $C\in \Sigma$ of positive measure and
$E(uf)=\lambda f$ for $0\neq \lambda \in \mathbb{C}$, then
$f=\frac{E(uf)}{\lambda}$, which means that $f$ is
$\mathcal{A}$-measurable. Therefore $E(uf)=E(u)f=\lambda f$ and
$(E(u)-\lambda)f=0$. This implies that $C\subseteq A_{\lambda}$
and so $\mu(A_{\lambda})>0$. This observations show that
$$\sigma_{p}(M_wEM_u)\setminus\{0\}=\{\lambda\in\mathbb{C}\setminus\{0\}:\mu(A_{\lambda,w})>0\},$$
where $A_{\lambda,w}=\{x\in X:E(uw)(x)=\lambda\}$.\\

 Here we show that a large class of integral operators are of the form of weighted conditional type operators. This means that we investigated centeredness and normality of integral operators on $L^2(\Sigma)$.\\

 Let $(X_1,\Sigma_1, \mu_1)$ and $(X_2,\Sigma_2, \mu_2)$ be two
$\sigma$-finite measure spaces and $X=X_1\times X_2$,
$\Sigma=\Sigma_1\times \Sigma_2$ and $\mu=\mu_1\times \mu_2$. Put
$\mathcal{A}=\{A\times X_2:A\in \Sigma_1\}$. Then $\mathcal{A}$ is
a sub-$\sigma$-algebra of $\Sigma$. Then for all $f$ in domain
$E^{\mathcal{A}}$ we have
$$E^{\mathcal{A}}(f)(x_1,x_2)=\int_{X_2}f(x_1,y)d\mu_2(y) \ \ \
\mu-a.e.$$ on $X$.\\

Also, if $(X,\Sigma, \mu)$ is a finite measure space and
$k:X\times X\rightarrow \mathbb{C}$ is a $\Sigma\otimes
\Sigma$-measurable function such that
$$\int_{X}|k(.,y)f(y)|d\mu(y)\in L^2(\Sigma)$$
for all $f\in L^2(\Sigma)$. Then the operator
$T:L^2(\Sigma)\rightarrow L^2(\Sigma)$ defined by
$$Tf(x)=\int_{X}k(x,y)f(y)d\mu, \ \ \ \ \ f\in L^2(\Sigma),$$
is called kernel operator on $L^2(\Sigma)$). We show that $T$ is a
weighted conditional type operator.\cite{gd} Since
$L^2(\Sigma)\times \{1\}\cong L^2(\Sigma)$ and
 $uf$ is a $\Sigma\otimes \Sigma$-measurable function, when $f\in
L^2(\Sigma)$. Then by taking $u:=k$ and $f'(x,y)=f(y)$, we get
that

$$E^{\mathcal{A}}(uf)(x,y)=E^{\mathcal{A}}(uf')(x,y)=\int_{X}u(x,y)f'(x,y)d\mu(y)=\int_{X}u(x,y)f(y)d\mu(y)=Tf(x).$$
Hence $T$ is a weighted conditional type operator.\\

By Theorem 2.1 we conclude that every well-defined integral
operator $T:L^2(\Sigma)\rightarrow L^0(\Sigma)$ such that
$$Tf(x)=\int_{X}k(x,y)f(y)d\mu, \ \ \ \ \ f\in L^2(\Sigma),$$
where $k:X\times X\rightarrow \mathbb{C}$ is a $\Sigma\otimes
\Sigma$-measurable function, is a bounded operator on
$L^2(\Sigma)$ if and only if
$$\int_{X}|k(.,y)|^2d\mu(y)<\infty, \mu \  a.e.$$\\

Now, we can derive conditions that
ensure the integral operators are centered on
$L^2(\Sigma)$-spaces as follows.\\

\vspace*{0.3cm} {\bf Theorem 3.1.} Let $Tf(x)=\int_{X}k(x,y)f(y)d\mu(y)$ for $f\in L^2(\Sigma)$. We have the followings:\\

(i) If $|\int_{X}k(x,y)d\mu(y)|^2=\int_{X}|k(x,y)|^2d\mu$ a.e. $\mu$, then the integral operator $T$ is centered on $L^2(\Sigma)$.\\

(ii) If the integral operator $T$ on $L^2(\Sigma)$ is centered, then $|\int_{X}k(x,y)d\mu(y)|^2=\int_{X}|k(x,y)|^2d\mu$ a.e. $\mu$ on $S(\int_{X}k(.,y)d\mu(y))$.\\

(iii) If $S(\int_{X}k(.,y)d\mu(y))=X$, then the integral operator $T$ on $L^2(\Sigma)$ is centered if and only if $|\int_{X}k(x,y)d\mu(y)|^2=\int_{X}|k(x,y)|^2d\mu$ a.e, $\mu$.\\

\vspace*{0.3cm} {\bf Theorem 3.2.} Let $Tf(x)=\int_{X}k(x,y)f(y)d\mu(y)$ for $f\in L^2(\Sigma)$. We have the followings:\\

(a) The integral operator $T$ is normal on $L^2(\Sigma)$ if and only if $$|\int_{X}k(x,y)d\mu(y)|^2=\int_{X}|k(x,y)|^2d\mu \ \ \ \  \  a.e. \mu.$$\\

(b) If $S(\int_{X}k(.,y)d\mu(y))=S(\int_{X}|k(.,y)|^2d\mu)$, the followings are equivalent.\\

(i) The integral operator $T$ on $L^2(\Sigma)$ is centered;\\

(ii) The integral operator $T$ on $L^2(\Sigma)$ is normal;\\

(iii) $|\int_{X}k(.,y)d\mu(y)|^2=\int_{X}|k(.,y)|^2d\mu$ a.e, $\mu$.\\

\vspace*{0.3cm} {\bf Theorem 3.3.} Let $T=M_wEM_u$ be bounded on $L^2(\Sigma)$ and $|E(uw)|^2=E(|u|^2)E(|w|^2)$ a.e, $\mu$. Then
 $$\sigma_{p}(M_wEM_u)=\sigma_{jp}(M_wEM_u).$$

 \vspace*{0.3cm} {\bf Proof.} Let $f\in L^2(\Sigma)\setminus \{0\}$ and $\lambda\in \mathbb{C}$, such that
 $wE(uf)=\lambda f$. Let $M=span\{f\}$, the closed linear subspace
 generated by $f$. Thus we can represent $T=M_wEM_u$ as the
 following $2\times 2$ operator matrix with respect to the
 decomposition $L^2(\Sigma)=M\oplus M^{\perp}$,
 $$
 T= \left[
         \begin{array}{rr}
              M_{\lambda} & PM_{w}EM_{u}- M_{\lambda} \\
              0 &  M_{w}EM_{u}-PM_{w}EM_{u}
          \end{array} \right].
$$
where $P$ is is the orthogonal projection of $L^2(\Sigma)$ onto
$M$. Since $|E(uw)|^2=E(|u|^2)E(|w|^2)$ a.e, $\mu$, for every $g\in
L^2(\Sigma)$ we have
$$\langle P|T^2|P(g)-P|T^{\ast}|^2P(g),g \rangle\geq0.$$

Hence $P(|T^2|-|T^{\ast}|^2)P\geq0$. Direct computation shows that
$$P|T^2|^2P= \left[
         \begin{array}{rr}
              M_{|\lambda|^4} & 0 \\
              0 &  0
          \end{array} \right]
.$$

Thus
$$\left[
         \begin{array}{rr}
              M_{|\lambda|^2} & 0 \\
              0 &  0
          \end{array} \right]\geq
         \left[
         \begin{array}{rr}
              M_{|\lambda|^2}+AA^{\ast} & 0 \\
              0 &  0
          \end{array} \right]
,$$ where $A=PM_{w}EM_{u}- M_{\lambda}$. This implies that $A=0$
and so $M$ is invariant under $T^{\ast}=M_{\bar{u}}EM_{\bar{w}}$,
too. Hence we conclude that
$T^{\ast}(f)=M_{\bar{u}}EM_{\bar{w}}(f)=\bar{\lambda} f$. This
means
that $\sigma_{p}(M_wEM_u)=\sigma_{jp}(M_wEM_u)$.\\

\vspace*{0.3cm} {\bf Remark 3.4.} Let $T=M_wEM_u$ be bounded on $L^2(\Sigma)$. Then;\\

(i) If $T$ is centered and $S(E(uw)E(u)E(w))=X$. Then  $\sigma_{p}(M_wEM_u)=\sigma_{jp}(M_wEM_u)$ and

$$\sigma_{jp}(M_wEM_u)\setminus \{0\}= \{\lambda\in\mathbb{C}\setminus\{0\}:\mu(A_{\lambda,w})>0\}.$$

(ii) If the operator $EM_u$ is centered and $S(E(u))=X$. We have $\sigma_{p}(EM_u)=\sigma_{jp}(EM_u)$ and

 $$\sigma_{jp}(EM_u)\setminus \{0\}= \{\lambda\in\mathbb{C}\setminus\{0\}:\mu(A_{\lambda})>0\}.$$

(iii) If the integral operator $Tf(x)=\int_{X}k(x,y)f(y)d\mu(y)$ is centered and $S(\int_{X}k(.,y)d\mu(y))=X$. We get that $\sigma_{p}(T)=\sigma_{jp}(T)$ and

$$\sigma_{jp}(T)\setminus \{0\}= \{\lambda\in\mathbb{C}\setminus\{0\}:\mu(A_{\lambda})>0\},$$

where $A_{\lambda}=\{x\in X:\int_{X}k(x,y)d\mu(y)=\lambda\}$.\\


\begin{thebibliography}{99}

\bibitem{afkz} J. Appell, E. V. Frovola, A. S. Kalitvin and P. P. Zabrejko,
Partial intergal operators on $C([a,b] \times [c,d])$,  Integral
Equations and Operator Theory {\bf 27} (1997), 125 - 140.

\bibitem{akn}  J. Appell, A.S. Kalitvin, and M.Z. Nashed, On some partial
integral equations arising in the mechanics of solids, Journal of
Applied Mathematics and Mechanics  {\bf 79} (1999), 703 - 713.

\bibitem{akz}  J. Appell, A.S. Kalitvin, and P.P. Zabrejko,  Partial integral
operators in Orlicz spaces with mixed norm,  Collo. Math. {\bf 78}
(1998), 293 - 306.

\bibitem{b} H. D. Brunk, On an extension of the concept conditional expectation, Proc. Amer. Math. Soc. {\bf 14} (1963), 298-304.


\bibitem{dhd}  P.G. Dodds, C.B. Huijsmans and B. De Pagter,
Characterizations of conditional expectation-type operators,
  Pacific J. Math. {\bf 141} (1990),
55-77.

\bibitem{dou}
 R. G. Douglas, Contractive projections on an $L\sb{1}$ space,
 Pacific J. Math. {\bf 15} (1965), 443-462.

\bibitem{e} Y. Estaremi, Essential norm of weighted conditional type operators on $L^p$-spaces, Positivity, to appear.


\bibitem{ej} Y. Estaremi and M.R. Jabbarzadeh, Weighted lambert type operators on
$L^{p}$-spaces, Oper. Matrices {\bf 1} (2013), 101-116.


\bibitem{gd} J.J. Grobler and B. de Pagter, Operators representable as
multiplication conditional expectation operators, Journal of
Operator Theory {\bf 48} (2002) 15 - 40.


\bibitem{her}
J. Herron, Weighted conditional expectation operators, Oper.
Matrices {\bf 1} (2011), 107-118.

\bibitem{iyy} M. Ito, T. Yamazaki and M. Yanagida, On the polar decomposition of the Aluthge transformation and related results, Journal of
Operator Theory {\bf 5} (2004) 303-319.

\bibitem{kz} A.S. Kalitvin and P.P. Zabrejko, On the theory of partial integral
operators, J. Integral Equations and Appl. {\bf 3} (1991), 351 -
382.

\bibitem{k} J.L. Kelley, Averaging operators on $C_{\infty}(X)$, Illinois J.
Math. {\bf 2} (1958), 214 - 223.

\bibitem{l} A. Lambert, A Hilbert $C^*$-module view of some spaces related
to probabilistic conditional expectation, Questiones Mathematicae
{\bf 22} (1999), 165 - 170.



\bibitem{mo}
Shu-Teh Chen, Moy,  Characterizations of conditional expectation
as a transformation on function spaces, Pacific J. Math. {\bf 4}
(1954), 47-63.


\bibitem{mm} B. B. Morrell and P. S. Muhly, centered operators.
Studia Math. {\bf 51} (1974), 251–263.



\bibitem{rao}
M. M. Rao, Conditional measure and applications, Marcel Dekker,
New York, 1993.

\bibitem{s} S. Sidak, On relations between strict  sense and wide sense
conditional expectation, Theory of Probability and Applications
{\bf 2} (1957), 267 - 271.


\bibitem{z}
A. C. Zaanen, Integration, 2nd ed., North-Holland, Amsterdam,
1967.



\end{thebibliography}
\end{document}